\documentclass[12pt]{article}
\usepackage{amsfonts,amssymb}
\usepackage{amsmath}
\usepackage{amsthm}
\newtheorem{theorem}{Theorem}
\newtheorem{lemma}{Lemma}
\newtheorem{prop}{Proposition}
\newtheorem*{proposition}{Proposition}
\newtheorem{corollary}{Corollary}

\title{On the Space of Symmetric Operators with Multiple Ground
States}
\author{A. Agrachev \thanks{SISSA, Trieste and MIAN, Moscow.}}
\date{}
\begin{document}
\maketitle

\begin{abstract} We study homological structure of the filtration
of the space of self-adjoint operators by the multiplicity of the
ground state. We consider only operators acting in a finite
dimensional complex or real Hilbert space but infinite dimensional
generalizations are easily guessed.
\end{abstract}

\section{Introduction}
This paper is dedicated to the memory of V. I. Arnold and is somehow
inspired by his works \cite{Ar}, \cite{Ar1} (see also \cite{ShVa}).
It opens a planned series of papers devoted to homological
invariants of the families of quadratic forms and related geometric
structures; Theorem~2 below forms a fundament of all further
development.

In this paper we study the filtration of the space of self-adjoint
operators by the multiplicity of the ground state. We restrict
ourself to operators in a finite dimensional complex or real
Hilbert space, but possible infinite dimensional generalizations
are easily guessed.

Let $\lambda_1(A)\le\lambda_2(A)\le\cdots\le\lambda_k(A)\le\cdots$
be the ordered eigenvalues of the operator $A$. The operators with
the ground state of multiplicity at least $k$ are characterized by
the equation $\lambda_1(A)=\lambda_k(A)$. Theorem~1 describes the
homotopy type of the space of nontrivial solutions to this
equation, it appears to be the Thom space of certain vector bundle
over a Grassmannian.

The growth of the multiplicity from $k$ to $k+1$ is realized by
the intersection of the space of solutions to the equation
$\lambda_1(A)=\lambda_k(A)$ with the space of solutions to the
equation $\lambda_k(A)=\lambda_{k+1}(A)$. We focus on the
homological structure of this intersection procedure.

As often happens, it is more convenient to accept the dual
viewpoint, i.\,e. to deal with the cohomology of the pair:
$\left(\, \mathbb B,\ \{A\in\mathbb B:
\lambda_1(A)\ne\lambda_k(A)\}\right)$ instead of the homology of
the space $\{A\in\mathbb B:\lambda_1(A)=\lambda_k(A)\}$, where
$\mathbb B$ is the space of all self-adjoint operators. Then the
intersections of cycles is substituted by the standard
cohomological product.

The space of solutions to the equation
$\lambda_k(A)=\lambda_{k+1}(A)$ is a cycle of codimension 3 in the
complex case and a cycle modulo 2 of codimension 2 in the real
case. The dual object is a 3-dimensional cohomology class in the
complex case and a 2-dimensional cohomology class modulo 2 in the
real one; we mean the class of the pair $\left(\,\mathbb B,\
\{A\in\mathbb B: \lambda_k(A)\ne\lambda_{k+1}(A)\}\right)$.

In both cases, we denote this cohomology class by $\Gamma_k$ and
study the map from the cohomology of the pair $\left(\,\mathbb B,\
\{A\in\mathbb B: \lambda_1(A)\ne\lambda_k(A)\}\right)$ to the
cohomology of the pair $\left(\,\mathbb B,\ \{A\in\mathbb B:
\lambda_1(A)\ne\lambda_{k+1}(A)\}\right)$ which sends any
cohomology class to its product with $\Gamma_k$. The main result
of the paper, Theorem~2 states that the sequence of these maps for
$k=1,2,\ldots$ is an exact sequence.

Let us consider the simplest case of self-adjoint operators
\linebreak on $\mathbb C^2$ or, in other words, of $2\times 2$
Hermitian matrices. The pair \linebreak $\left(\,\mathbb B,\
\{A\in\mathbb B: \lambda_1(A)\ne\lambda_k(A)\}\right)$ equals
$\left(\mathbb R^4,\emptyset\right)$ for $k=1$ and $\left(\mathbb
R^4,\mathbb R^4\setminus\mathbb R \right)$ for $k=2$. The exact
sequence of Theorem~2 is reduced to the obvious sequence
$$
0\to H^*(\mathbb R^3)\to H^{*+3}(\mathbb R^3,\mathbb R^3\setminus
0)\to 0.
$$

The general multidimensional calculation is far from being trivial
and has perhaps a fundamental nature as we hope to show in the
forthcoming publications.

In the next section we introduce some notations and recall
well-known facts on the spaces of self-adjoint operators. The
Theorems 1 and 2 are formulated and proved in Section~3.

\smallskip
All pairs of topological spaces and their subspaces we deal with
are homotopy equivalent to pairs of finite cell complexes and
their subcomplexes. Let $(M,X),\ (M,Y),\ (M,X\cup Y)$ be such
pairs, $\xi\in H^i(M,X),\ \eta\in H^j(M,Y)$; then
$\xi\smile\eta\in H^{i+j}(M,X\cup Y)$ is the cohomological product
of $\xi$ and $\eta$.

\section{Preliminaries}

We consider the spaces of self-adjoint operators on $\mathbb R^n$ or
$\mathbb C^n$. In both cases, given an operator $A$ we denote by
$\lambda_1(A)\le\cdots\le\lambda_n(A)$ its ordered eigenvalues. Let
$I$ be the unit operator and $\alpha$ a positive real number.
Obviously, $\lambda_i(A\pm\alpha I)=\lambda_i(A)\pm\alpha,\
\lambda_i(\alpha A)=\alpha\lambda_i(A)$. Moreover, $A\pm\alpha I$
and $\alpha A$ have the same eigenvectors as $A$. It is convenient
do not distinguish the operators obtained one from another by just
described trivial transformations.

We denote by $\mathbb S(\mathbb R)$ (correspondingly by $\mathbb
S(\mathbb C)$) the space of all non-scalar self-adjoint linear
operators $A:\mathbb R^n\to\mathbb R^n$ (correspondingly
$A:\mathbb C^n\to\mathbb C^n)$ factorized by the equivalence
relation $A\sim(\alpha A+\beta I),\
\forall\,\alpha>0,\beta\in\mathbb R$. Then $\mathbb S(\mathbb R)$
is homeomorphic to the sphere $S^{\frac {n(n+1)}2-2}$ and $\mathbb
S(\mathbb C)$ is homeomorphic to $S^{n^2-2}$. In what follows we
often deal simultaneously with the real and Hermitian cases and
simply omit the argument of $\mathbb S$. The join of $\mathbb S$
and the origin is denoted by $\mathbb B$; this is the ball of
dimension $\frac {n(n+1)}2-1$ in the real case and the ball of
dimension $n^2-1$ in the Hermitian case.

We can substitute the factorization by the normalization and
define $\mathbb S$ as the space of self-adjoint operators $A$ such
that $\sum\limits_{i=1}^n\lambda_i(A)=0,\
\sum\limits_{i=1}^n\lambda_i^2(A)=1$; then $\mathbb B$ is defined
by the relations $\sum\limits_{i=1}^n\lambda_i(A)=0,\
\sum\limits_{i=1}^n\lambda_i^2(A)\le 1$. Anyway, the normalization
is sometimes less convenient than the factorization and we often
use the same symbols for the equivalence classes and their
representatives; this simplifies notations and does not lead to a
confusion.

Now consider open subsets
$$
\Sigma_{k,k+1}\doteq\{A\in\mathbb
S:\lambda_k(A)\ne\lambda_{k+1}(A)\}.
$$
The following facts are well-known:

\begin{prop} $\mathbb S\setminus\Sigma_{k,k+1}$ is an algebraic
subset of $\mathbb S$ of codimension 2 in the real case and of
codimension 3 in the Hermitian case. Singular locus of $\mathbb
S\setminus\Sigma_{k,k+1}$ consists of the operators with at least
triple eigenvalue $\lambda_k$; it is an algebraic subset of
$\mathbb S$ of codimension 5 in the real case and of codimension 8
in the Hermitian case. Moreover, regular part of $\mathbb
S\setminus\Sigma_{k,k+1}$ is orientable in the Hermitian case.
\end{prop}
{\bf Sketch of the proof.} Let $A_0\in\mathbb
S\setminus\Sigma_{k,k+1}$ and $J_{A_0}$ be the set of all
\linebreak $j\in\{1,\ldots,n\}$ such that
$\lambda_j(A_0)=\lambda_k(A_0)$; then $\#J_{A_0}\ge 2$. Given a
self-adjoint operator $A$, we set $ K_A=span\{x\in X:
Ax=\lambda_j(A)x,\ j\in J_{A_0}\},$ where $X$ is $\mathbb R^n$ in
the real case and $X$ is $\mathbb C^n$ in the Hermitian case.

Let $\mathcal O_0$ be a neighborhood of $A_0$ in $\mathbb S$ such
that $K_A\cap K_{A_0}^\perp=0,\ \forall\,A\in\mathcal O_0$. We
denote by $P^k_{A_0A}:K_A\to K_{A_0}$ the restriction to $K_A$ of
the orthogonal projector $X\to K_{A_0}$ and set
$\Phi(A)=P_{A_0A}AP^{-1}_{A_0A},\ A\in\mathcal O_0$. Then $\Phi$
is a well-defined rational map from $\mathcal O_0$ to the space of
self-adjoint operators on $K_{A_0}$. The differential of $\Phi$ at
$A_0$ sends $A$ to the composition of $A\bigr|_{K_{A_0}}$ with the
orthogonal projection $X\to K_{A_0}$. In particular, $D_{A_0}\Phi$
is surjective; hence $\Phi$ is a submersion on a neighborhood of
$A_0$. We may assume that $\Phi$ is a submersion on the whole
$\mathcal O_0$. Moreover, $\lambda_i(\Phi(A))=\lambda_{i+j_0}(A),\
i=1,\ldots,\#J_{A_0}$, where $j_0=\min J_{A_0}$.

Let $A\in\mathcal O_0$; we obtain that $J_A=J_{A_0}$ if and only
if $\Phi(A)$ is a scalar operator. Hence $\{A\in\mathcal O_0:
J_A=J_{A_0}\}$ is a regular algebraic subset of $\mathcal O_0$ of
codimension $\frac{j_0(j_0+1)}2-1$ in the real case and $j_0^2-1$
in the Hermitian case.

It remains to prove the orientability in the Hermitian case. It is
sufficient to show that the space of self-adjoint operators on
$K_A,\ A\in\mathbb S$ has a canonical orientation. The orientation
of the space of self-adjoint operators is induced by the orientation
of the space $K_A$ itself, and the orientation of $K_A\subset\mathbb
C^n$ is defined by the complex structure (any complex space has a
canonical orientation). \qquad $\square$

\medskip
Proposition 1 implies that $H_{\dim\mathbb S-2}(\mathbb
S\setminus\Sigma_{k,k+1};\mathbb Z_2)=\mathbb Z_2$ in the real
case and $H_{\dim\mathbb S-3}(\mathbb
S\setminus\Sigma_{k,k+1};\mathbb Z)=\mathbb Z$ in the Hermitian
case. According to the Alexander duality,
$H^1(\Sigma_{k,k+1};\mathbb Z_2)=\mathbb Z_2$ in the real case and
the generator of $H^1(\Sigma_{k,k+1};\mathbb Z_2)$ applied to a
closed curve in $\Sigma_{k,k+1}$ equals the linking number of the
curve and $\mathbb S\setminus\Sigma_{k,k+1}$ modulo 2. Similarly,
$H^2(\Sigma_{k,k+1};\mathbb Z)=\mathbb Z$ in the Hermitian case
and the generator of $H^1(\Sigma_{k,k+1};\mathbb Z)$ applied to a
compact oriented surface in $\Sigma_{k,k+1}$ equals the linking
number of the surface and \linebreak $\mathbb
S\setminus\Sigma_{k,k+1}$.

Further in this paper we always consider homology and cohomology
with coefficients in $\mathbb Z_2$ in the real case and with
coefficients in $\mathbb Z$ in the Hermitian case, and we omit the
indication of coefficients in order to simplify notations. We also
denote by $\varepsilon$ the codimension of $\mathbb
S\setminus\Sigma_{k,k+1}$ in $\mathbb S$; then $\varepsilon=2$ in
the real case and $\varepsilon=3$ in the Hermitian case.

Given $A\in\Sigma_{k,k+1}$, we set
$$
E^k_A=span\{x\in X: Ax=\lambda_i(A)x,\ i=1,\ldots,k\},
$$
where $X=\mathbb R^n$ in the real case and $X=\mathbb C^n$ in the
complex case. Then $\mathcal E^k=\{(x,A): A\in\Sigma_{k,k+1}, x\in
E^k_A\}$ is a $k$-dimensional vector subbundle of the trivial
bundle $X\times\Sigma_{k,k+1}$ over $\Sigma_{k,k+1}$. Let
$\gamma_k\in H^{\varepsilon-1}(\Sigma_{k,k+1})$ be the first
Stiefel--Whitney characteristic class of this bundle in the real
case and the first Chern characteristic class in the Hermitian
case.

\begin{prop} $\gamma_k$ is a generator of
$H^{\varepsilon-1}(\Sigma_{k,k+1})$.
\end{prop}
{\bf Proof.} We have to compute the characteristic classes of the
restriction of the bundle $\mathcal E^k\to\Sigma_{k,k+1}$ to a
$(\varepsilon-1)$-dimensional compact submanifold of
$\Sigma_{k,k+1}$ whose linking number with $\mathbb
S\setminus\Sigma_{k,k+1}$ equals $\pm 1$.

We denote by $\mathbb S_2$ the space of self-adjoint operators $B$
on $\mathbb R^2$ (in the real case) or on $\mathbb C^2$ (in the
Hermitian case) such that $\lambda_1(B)+\lambda_2(B)=0,\
\lambda^2_1(B)+\lambda^2_2(B)=1$ and set $\mathbb B_2=conv(\mathbb
S_2)$. Then $\mathbb S_2$ is a $(\varepsilon-1)$-dimensional
sphere and $\mathbb B_2$ is a $\varepsilon$-dimensional ball,
$\mathbb S_2=\partial\mathbb B_2$. Let $A_{-}$ be a self-adjoint
operator on $\mathbb R^{k-1}$ (or on $\mathbb C^{k-1}$) with
simple eigenvalues such that $\lambda_{k-1}(A_{-})<-1$ and $A_+$
be a self-adjoint operator on $\mathbb R^{n-k-1}$ (or on $\mathbb
C^{n-k-1}$) with simple eigenvalues such that $\lambda_1(A_+)>1$.
Then $A_{-}\oplus\mathbb S_2\oplus A_+$ is the required
$(\varepsilon-1)$-dimensional submanifold of $\mathbb S$. Indeed,
$A_{-}\oplus\mathbb S_2\oplus A_+=\partial\left(A_{-}\oplus\mathbb
B_2\oplus A_+\right)$ and $\left(A_{-}\oplus\mathbb B_2\oplus
A_+\right)\cap\left(\mathbb
S\setminus\Sigma_{k,k+1}\right)=\left(A_{-}\oplus 0\oplus
A_+\right)$; moreover, the intersection is transversal. Hence the
linking number of $A_{-}\oplus\mathbb S_2\oplus A_+$ and $\mathbb
S\setminus\Sigma_{k,k+1}$ equals $\pm 1$.

The restriction of the bundle $\mathcal E^k\to\Sigma_{k,k+1}$ to
$A_{-}\oplus\mathbb S_2\oplus A_+$ splits in the sum of a trivial
vector bundle and a linear bundle over $\mathbb S_2$ whose fiber
at $B\in\mathbb S_2$ is the eigenspace of the eigenvalue
$\lambda_1(B)$. It is easy to see that the map sending
$B\in\mathbb S_2$ to the eigenspace of the eigenvalue
$\lambda_1(B)$ is the diffeomorphism of $\mathbb S_2$ and the
projective line (real or complex). This diffeomorphism identifies
our linear bundle with the tautological bundle of the projective
line. \qquad $\square$

\medskip
Let $\delta:H^i(\Sigma_{k,k+1})\to H^{i+1}(\mathbb
B,\Sigma_{k,k+1})$ be the isomorphism induced by the exact
cohomological sequence of the pair $(\mathbb B,\Sigma_{k,k+1})$. We
set $\Gamma_k=\delta\circ\gamma_k\in H^\varepsilon(\mathbb
B,\Sigma_{k,k+1})$. The value of $\Gamma_k$ on a relative cycle
$\xi\in C_\varepsilon(\mathbb B,\Sigma_{k,k+1})$ is the intersection
number of $\xi$ and $conv(\mathbb S\setminus\Sigma_{k,k+1})$.

\medskip
We conclude this section with an explicit expression for a closed
two-form representing the class $\gamma_k$ in the Hermitian case.

Let $A$ be a self-adjoint operator with simple eigenvalues and
$e_1,\ldots,e_n$ an orthonormal basis of its eigenvectors such that
$Ae_i=\lambda_i(A)e_i,\ i=1,\ldots,n$. Then $e_i$ are defined up to
a complex multiplier of the absolute value 1. Let
$\langle\cdot,\cdot\rangle$ denotes the Hermitian product and $B$ be
another self-adjoint operator. It is easy to see that the wedge
square over $\mathbb R$ of the complex number $\langle
Be_i,e_j\rangle$ depends only on $A, B, i, j$ and not on the choice
of the eigenvectors. In particular, $\bigwedge^2_{\mathbb R}\langle
dA\,e_i,e_j\rangle$ is a well-defined two-form on the space of
self-adjoint operators with simple eigenvalues. We have: $
\bigwedge^2_{\mathbb R}\langle
dA\,e_i,e_j\rangle\left(\frac\partial{\partial
B_1},\frac\partial{\partial B_2}\right)=\det_{\mathbb
R}\left(\langle B_1e_i,e_j\rangle,\langle B_2e_i,e_j\rangle\right),
$ where complex numbers are treated as vectors in $\mathbb R^2$.

\begin{prop} The form
$\Omega_k=\sum\limits_{i=1}^k\sum\limits_{j=k+1}^n\frac
1{\pi(\lambda_i-\lambda_j)^2}\bigwedge^2_{\mathbb R}\langle
dA\,e_i,e_j\rangle$ represents the restriction of $\gamma_k$ to the
space of self-adjoint operators with simple eigenvalues\footnote{The
form $\Omega_k$ is locally bounded in the topology of
$\Sigma_{k,k+1}$. Moreover, any two-dimensional cycle in
$\Sigma_{k,k+1}$ is homotopic to a cycle in the space of
self-adjoint operators with simple eigenvalues. Hence the form
$\Omega_k$ indeed represents $\gamma_k$.}.
\end{prop}
{\bf Sketch of the proof.} We have to demonstrate that $\Omega_k$
represents the first Chern class of the vector bundle $\mathcal
E^k$ restricted to the space of self-adjoint operators with simple
eigenvalues. This restriction of $\mathcal E^k$ is the direct sum
of line bundles $\mathcal L^i,\ i=1,\ldots,k$, where the fiber of
$\mathcal L^i$ at $A$ is the line \linebreak
$L^i_A\doteq\{z\in\mathbb C^n:Az=\lambda_i(A)z\}$. We have to show
that $\Omega_k$ represents the class
$\sum\limits_{i=1}^kc_1(\mathcal L^i)$.

Consider the associated with $\mathcal L^i$ principal $S^1$-bundle
$\mathcal C^i$ whose fiber at $A$ is $C^i_A\doteq\{e_i\in\mathbb
C^n:e_i\in L^i_A,\ |e_i|=1\}$. Given a smooth curve $t\mapsto
A(t)$ in the space of self-adjoint operators with simple
eigenvalues and $e_i(0)\in C^i_{A(0)}$, the condition $\langle
\dot e_i(t),e_i(t)\rangle=0$ defines a canonical lift $t\mapsto
e_i(t)$ of the curve $A(\cdot)$ to the bundle $\mathcal C^i$.
These lifts are the parallel translations for a connection on the
bundle $\mathcal C^i$ along curves in the base space of the
bundle. The form of this connection equals $\Im\langle
de_i,e_i\rangle$, where $\Im$ denotes the imaginary part of a
complex number.

The exterior differential of the form $\Im\langle de_i,e_i\rangle$
is the pullback of the curvature form $R_i$ of the connection. An
immediate calculation shows that $R_i\left(\frac\partial{\partial
B_1},\frac\partial{\partial
B_2}\right)\Bigr|_A=2\Im\left\langle\frac{\partial e_i(A)}{\partial
B_2},\frac{\partial e_i(A)}{\partial B_1}\right\rangle$, where
$\frac{\partial e_i(A)}{\partial B}=\frac
d{dt}e_i(A+tB)\bigr|_{t=0}$ and $t\mapsto e_i(A+tB)$ is parallel
along the curve $t\mapsto A+tB$. The differentiation by $t$ of the
equation $
\left\langle(A+tB)e_i(A+tB),e_j(A)\right\rangle=\lambda_i(A+tB)\langle
e_i(A+tB),e_j(A)\rangle$ gives: $\left\langle\frac{\partial
e_i(A)}{\partial B}, e_j(A)\right\rangle=\frac
1{\lambda_i(A)-\lambda_j(A)}\langle Be_i(a),e_j(A)\rangle,\ \forall
j\ne i.$ Hence $\frac{\partial e_i}{\partial B}=\sum\limits_{i\ne
j}\frac 1{\lambda_i-\lambda_j}\langle Be_i,e_j\rangle e_j$ and
$R_i\left(\frac\partial{\partial B_1},\frac\partial{\partial
B_2}\right)=\sum\limits_{j\ne i}\frac
2{(\lambda_i-\lambda_j)^2}\Im\left(\langle B_2e_i,e_j\rangle\langle
B_1e_i,e_j\rangle\right)$.

On the other hand $\Im(z^1z^2)=\det_{\mathbb R}(z^2,z^1)$ for any
complex numbers $z^1,z^2$. We obtain: $R_i=\sum\limits_{j\ne
i}\frac 2{(\lambda_i-\lambda_j)^2}\bigwedge_{\mathbb R}^2\langle
dA\,e_i,e_j\rangle$. Summing up, we get the desired expression for
the form $\Omega_k=\frac 1{2\pi}\sum\limits_{i=1}^kR_i$
representing the class $\gamma_k. \qquad \square $

\section{Main Results}

We are going to study the filtrations
$$
M_k=\{A\in\mathbb S: \lambda_1(A)=\lambda_{k+1}(A)\},\ M^k=\mathbb
S\setminus M_k, \quad k=0,\ldots,n-1.
$$
It is easy to see that
$$
W_k\doteq\{A\in\mathbb S: \lambda_{k+1}(A)=\lambda_n(A)\}\subset
M^k
$$
is deformation retract of $M^k$. The retraction $\phi_k:M^k\to
W_k$ changes only eigenvalues of the operators while the
eigenvectors are kept fixed.

The involution $A\mapsto (-A),\ A\in\mathbb S$, transforms $W_k$
into $M_{n-k-1}$; hence $M^k$ is homotopy equivalent to
$M_{n-k-1}$. Note also that the map $A\mapsto A-\lambda_1(A)I,\
A\in\mathbb S$, induces homeomorphism of $M_k$ and the space of
nonzero nonnegative self-adjoint operators of rank $<n-k$
factorized by the equivalence relation $A\sim\alpha A,\
\forall\alpha>0$.

In what follows, $Gr_k(m)$ is the Grassmannian of $k$-dimensional
subspaces of $\mathbb R^m$ or $\mathbb C^m$.
\begin{theorem} $M^k$ has homotopy type of the Thom space of a
real vector bundle over the Grassmannian $Gr_k(n-1)$; the
dimension of the bundle equals $\frac{k(k+1)}2+k-1$ in the real
case and $k^2+2k-1$ in the Hermitian case.
\end{theorem}
{\bf Proof.} Let $e$ be a unit length vector (in $\mathbb R^n$ or
$\mathbb C^n$). We set
$$
G_k(e)=\{A\in M_{k-1}\cap W_k : Ae=\lambda_n(A)e\}.
$$
Then $G_k(e)\cong Gr_k(n-1)$. Moreover, a neighborhood of $G_k(e)$
in $W_k$ is a smooth manifold and the normal bundle of $G_k(e)$ in
this manifold has dimension $\frac{k(k+1)}2+k-1$ in the real case
and dimension $k^2+2k-1$ in the Hermitian case. In fact, the
normal bundle splits in the sum of two subbundles. The first one
is the normal bundle of $G_k(e)$ in $\{A\in W_k:
Ae=\lambda_n(A)e\}$; it is isomorphic to the bundle of the
self-adjoint endomorphisms with zero trace of the tautological
bundle of $Gr_k(n-1)$. The second one is the normal bundle of
$G_k(e)$ in $M_{k-1}\cap W_k$; it is isomorphic to the normal
subbundle of $Gr_k(n-1)$ in $Gr_k(n)$. The theorem is an immediate
corollary of the following
\begin{lemma} $W_k\setminus G_k(e)$ is a contractible space.
\end{lemma}
{\bf Proof.} We'll contract $W_k\setminus G_k(e)$ to the point
$-e^*\otimes e\in W_k\setminus G_k(e)$. The contraction sends
$(A,t)\in\left(W_k\setminus G_k(e)\right)\times[0,1]$ to
$\phi_k(A_t)$, where $A_t=(1-t)A-te^*\otimes e$.

It remains to prove that the contraction is correctly defined,
i.\,e. that $A_t\in M^k$ and $\phi_k(A_t)$ does not belong to
$G_k(e)$, $\forall t\in[0,1]$. We have: $\phi_k(A_t)\in G_k(e)$ if
and only if $A_t\in M_{k-1}$ and $e$ is orthogonal to $\{x\in X:
A_tx=\lambda_1(A_t)x\}$, where $X$ is $\mathbb R^n$ (in the real
case) or $\mathbb C^n$ (in the Hermitian case).

Set $A^t=A-\frac t{1-t}e^*\otimes e=\frac 1{1-t}A_t$; the positive
multiplier does not influence the multiplicity of the eigenvalues
and we will work with $A^t$ instead of $A_t$. We consider
separately two cases.

\medskip\noindent
{\bf 1.} $e\perp \{x\in X: Ax=\lambda_i(A)x,\ i=1,\ldots k\}$.
Then $e$ is an eigenvector of $A,\ Ae=\lambda_n(A)e$. Hence all
$A^t$ have common eigenvectors. Moreover, $(n-1)$ eigenvalues of
$A^t$ (counted with the multiplicities) are equal to eigenvalues
of $A$ while the eigenvalue corresponding to the eigenvector $e$
is monotone decreasing from $\lambda_n(A)$ to $-\infty$ as $t$
runs from 0 to 1. Besides that, $\lambda_1(A)\ne\lambda_k(A)$
since $A\notin G_k(e)$; hence $A^t\in M^k$.

The equality $\lambda_1(A^t)=\lambda_k(A^t)$ is valid for some
$t\in(0,1]$ if and only if $\lambda_1(A)=\lambda_{k-1}(A)$; then
$A^te=\lambda_1(A^1)e$ and $\phi_k(A^t)\notin G_k(e)$.

\medskip\noindent
{\bf 2.} $e\not\perp \{x\in X: Ax=\lambda_i(A)x,\ i=1,\ldots k\}$.
The restriction of the quadratic form $x\mapsto\langle
A^tx,x\rangle$ to the hyperplane $e^\perp$ does not depend on $t$.
Hence the minimax principle for the eigenvalues implies:
$$
\lambda_1(A^t)\le\lambda_1(A)\le\lambda_2(A^t)\le\cdots\le\lambda_k(A)\le\lambda_{k+1}(A^t).
$$
Assume that $\lambda_1(A^t)=\lambda_{k+1}(A^t)$. Hence
$$
\lambda_1(A)=\lambda_k(A)=\lambda_1(A^t)=\min\limits_{|x|=1}\langle
A^tx,x\rangle.
$$
At the same time,
$$
\langle A^tx,x\rangle=\lambda_1(A)-\frac t{1-t}\langle
e,x\rangle^2,\quad \forall x\in\{x\in X: Ax=\lambda_1(A)x\}.
$$
We obtain the contradiction with the assumption
$$
e\not\perp\{x\in X:Ax=\lambda_1(A)x\}=\{x\in X:Ax=\lambda_i(A)x,\
i=1,\ldots,k\}.
$$
Hence $A_t\in M^k$.

Now assume that $\lambda_1(A^t)=\lambda_k(A^t)$ and $e\perp\{x\in
X:A^tx=\lambda_1(A^t)x\}$. Then
$$
\lambda_1(A^t)|x|^2=\langle A^tx,x\rangle=\langle
Ax,x\rangle=\lambda_1(A)|x|^2, $$ $$ \forall x\in\{x\in
X:A^tx=\lambda_i(A^t)x,\ i=1,\ldots,k\}.
$$
Hence $\lambda_1(A)=\lambda_k(A)$ and
$$
\{x\in
X:A^tx=\lambda_i(A^t)x,\ i=1,\ldots,k\}=
$$
$$
\{x\in X:Ax=\lambda_i(A)x,\ i=1,\ldots,k\}.
$$
We obtain the
contradiction with the assumption {\bf 2.} \qquad $\square$

\medskip
Let $u_k\in H^{\nu_k}(M^k)$ be the Thom class of the normal bundle
of $G_k(e)$ in $W_k,\ \nu_k=\frac{k(k+1)}2+k-1$ in the real case
and $\nu_k=k^2+2k-1$ in the Hermitian case. Let $\mathcal G_k$ be
the total space of this bundle and $G_k(e)\subset\mathcal G_k$ its
zero section. We have:
$$
\tilde H^\cdot(M^k)=H^\cdot(\mathcal G_k,\mathcal G_k\setminus
G_k(e)), \quad H^\cdot(G_k(e))=H^\cdot(\mathcal G_k)
$$
and the cohomology product of the classes from $\tilde
H^\cdot(M^k)$ and $H^\cdot(G_k(e))$ is a well-defined class in
$\tilde H^\cdot(M^k)$. Then $\xi\mapsto u_k\smile\xi,\ \xi\in
H^\cdot(G_k(e)),$ is the Thom isomorphism of $H(G_k(e))$ and
$\tilde H(M^k)$. Recall that $u_k\bigr|_{G_k(e)}\in
H^{\nu_k}(G_k(e))$ is the Euler class of the bundle $\mathcal
G_k\to G_k(e)$.

\begin{lemma} $u_k\bigr|_{G_k(e)}=0$.
\end{lemma}
{\bf Proof.} The bundle $\mathcal G_k$ splits in the sum of two
subbundles as it was explained in the proof of Theorem~1. We'll
prove that the first subbundle, i.\,e. the bundle of self-adjoint
endomorphisms with zero trace of the tautological bundle of the
Grassmannian has zero Euler class. It is sufficient to show that
the induced bundle over the flag space has zero Euler class. This
is easy. Indeed, the bundle over the flag space has natural
non-vanishing sections: the value of such a section at a flag is
the self-adjoint operator with prescribed simple eigenvalues whose
eigenspaces are the elements of the flag.\qquad $\square$

\begin{corollary} The cohomology product of any two elements
of $\tilde H^\cdot(M^k)$ is zero.
\end{corollary}
{\bf Proof.} Due to the Thom isomorphism, it is sufficient to show
that \linebreak $u_k\smile u_k=0$, but $u_k\smile u_k$ is the
image of $u_k\bigr|_{G_k(e)}=0$ under the Thom isomorphism.\qquad
$\square$

\medskip
Obviously, $M^k=M^{k-1}\cup\Sigma_{k,k+1}$. We consider the
homomorphisms
$$
\mathbf d_k:H^\cdot(\mathbb B,M^{k-1})\to H^\cdot(\mathbb
B,M^k),\quad k=1,\ldots,n-1,
$$
acting by multiplication with the class $\Gamma_k\in
H^{\varepsilon}(\mathbb B,\Sigma_{k,k+1})$ defined at the end of
Section~2:
$$
\mathbf d_k(\xi)=\Gamma_k\smile\xi,\quad \xi\in H^\cdot(\mathbb
B,M^{k-1}).
$$
Recall that $\varepsilon=2$ in the real case and $\varepsilon=3$
in the Hermitian case.

\begin{theorem}
$$
0\to H^\cdot(\mathbb B)\stackrel{\mathbf d_1}{\to}H^\cdot(\mathbb
B,M^1)\stackrel{\mathbf d_2}{\to}\cdots\stackrel{\mathbf
d_{n-2}}{\to}H^\cdot(\mathbb B,M^{n-2})\stackrel{\mathbf
d_{n-1}}{\to}H^\cdot(\mathbb B,\mathbb S)\to 0 \eqno (1)
$$
is an exact sequence.
\end{theorem}

{\bf Proof.} We make calculations only for the real case; the
Hermitian version is obtained by obvious modifications.

First note that $\Gamma_k\smile\Gamma_{k-1}=0$. Indeed,
$\Gamma_k\smile\Gamma_{k-1}$ is an element of $H^4(\mathbb
B,\Sigma_{k-1,k}\cup\Sigma_{k,k+1})=H^4(\mathbb
S,\Sigma_{k-1,k}\cup\Sigma_{k,k+1})$ but
$$
\mathcal
S\setminus(\Sigma_{k-1,k}\cup\Sigma_{k,k+1})=\{A\in\mathcal S:
\lambda_{k-1}(A)=\lambda_{k+1}(A)\}
$$
is a codimension 5 algebraic subset of $\mathcal S$ (see
Proposition~1). Hence $\mathbf d_k\circ\mathbf d_{k-1}=0$ and (1)
is a cochain complex. We have to prove that this complex has
trivial cohomology.

Consider the spaces:
$$
\Sigma_{1,k,k+1}\doteq\{A\in\mathbb
S:\lambda_{1}(A)\ne\lambda_k(A)\ne\lambda_{k+1}(A)\}=M^{k-1}\cap\Sigma_{k,k+1}.
$$
Then $\phi_k\left(\Sigma_{k,k+1}\right)=\Sigma_{k,k+1}\cap W_k$
and $\Sigma_{k,k+1}\cap W_k$ is a deformation retract of
$\Sigma_{k,k+1}$. Similarly,
$\phi_k\left(\Sigma_{1,k,k+1}\right)=\Sigma_{1,k,k+1}\cap W_k$ and
$\Sigma_{1,k,k+1}\cap W_k$ is a deformation retract of
$\Sigma_{1,k,k+1}$.

The map $A\mapsto span\{x\in\mathbb R^n: Ax=\lambda_i(A)x,\
i=1,\ldots,k\}$ endows the space $\Sigma_{k,k+1}\cap W_k$ with the
structure of the fiber bundle over $Gr_k(n)$, where the fiber at
$E\in Gr_k(n)$ is the space of all self-adjoint operators $A:E\to
E$ such that $\sum\limits_{i=1}^k\lambda_i(A)\le 0,\
\sum\limits_{i=1}^k\lambda^2_i(A)+ \frac 1{n-k}
\left(\sum\limits_{i=1}^k\lambda_i(A)\right)^2=1$: such operators
are uniquely extended to (normalized) operators from $W_k$. The
fiber is thus a ball of dimension $\frac{k(k+1)}2 -1=\nu_{k-1}+1$.

Moreover, $\left(\Sigma_{k.k+1}\setminus
\Sigma_{1,k,k+1}\right)\cap W_k$ is a section of the bundle
$\Sigma_{k,k+1}\cap W_k\to Gr_k(n)$, where the value of the
section at $E\in Gr_k(n)$ is a scalar operator on $E$, ``the
center of the ball". Hence ``the spherical bundle"
$$
\{A\in \Sigma_{k,k+1}\cap W_k:
\sum\limits_{i=1}^k\lambda_i(A)=0\}
$$
with a typical fiber
$S^{\nu_{k-1}}$ is a homotopy retract of $\Sigma_{1,k,k+1}$. Let
\begin{multline*} S_E^{\nu_{k-1}}\doteq\{A\in
\Sigma_{k,k+1}\cap W_k: \sum_{i=1}^k\lambda_i(A)=0,\ Ax_i=\lambda_i(A)x_i,\\
x_i\in E\setminus\{0\}, i=1,\ldots, k\}\quad (2)
\end{multline*}
be the fiber at $E$ of this spherical bundle.

\begin{lemma} The restriction $u_{k-1}\bigr|_{S_E^{\nu_{k-1}}}$ of
the class $u_{k-1}\in H^{\nu_{k-1}}(M^{k-1})$ induced by the
inclusion
$$
M^{k-1}\supset \Sigma_{1,k,k+1}\supset S_E^{\nu_{k-1}}
$$
is the generator of $H^{\nu_{k-1}}(S^{\nu_{k-1}}_E)$.
\end{lemma}
{\bf Proof.} The value of the Thom class $u_{k-1}$ on the cycle
$S^{\nu_{k-1}}_E$ is the intersection number of the cycle
$\phi_{k-1}(S^{\nu_{k-1}}_E)$ with $G_{k-1}(e)$ in $M^{k-1}\cap
W_{k-1}$. Obviously, this number does not depend on $E$. Take $E$
such that $e\not\perp E$. Then the intersection of
$\phi_{k-1}(S^{\nu_{k-1}}_E)$ and $G_{k-1}(e)$ is transversal and
consists of one point $A_0$ characterized by the relations
$$
A_0\in G_{k-1}(e),\quad \{x\in\mathbb R^n:
A_0x=\lambda_1(A_0)x\}=e^\perp\cap E. \eqno \square
$$

\begin{corollary} Let $v_{k-1}=u_{k-1}\bigr|_{\Sigma_{1,k,k+1}}$.
Then the ring $H^\cdot(\Sigma_{1,k,k+1})$ is a free module over
the ring $H^\cdot(Gr_k(n))$ with the basis $1,\, v_{k-1}$.
Moreover, $v_{k-1}\smile v_{k-1}=0$.
\end{corollary}
{\bf Proof.} The module structure is induced by the bundle
structure \linebreak $\Sigma_{1,k,k+1}\cap W_k\to Gr_k(n).$ The
fact that the module is free follows from Lemma~3 and the
Leray--Hirsch theorem. The equality $v_{k-1}\smile v_{k-1}=0$
follows from the equality $u_{k-1}\smile u_{k-1}=0$ (see
corollary~1). \qquad $\square$

\begin{lemma} The inclusions $M^{k-1}\subset M^k$ and $\Sigma_{k,k+1}\subset M^k$
induce zero homomorphisms of the reduced cohomology groups.
\end{lemma}
{\bf Proof.} We have: $\phi_k(M^{k-1})\subset M^k\setminus
G_k(e)$. Hence $M^{k-1}$ is contained in the contractible subset
of $M^k$ and the restriction to $M^{k-1}$ makes trivial any
cohomology class from $\tilde H(M^k)$. Now consider the inclusions
$$
\Sigma_{1,k,k+1}\subset \Sigma_{k,k+1}\subset M^k.
$$
Corollary 2 implies that the inclusion $\Sigma_{1,k,k+1}\subset
\Sigma_{k,k+1}\cong Gr_k(n)$ induces the injective homomorphism
$H^\cdot(\Sigma_{k,k+1})\to H^\cdot(\Sigma_{1,k,k+1})$. On the
other hand, $\phi_k(\Sigma_{1,k,k+1})\subset M^k\setminus G_k(e)$.
Hence the composition of the induced by the inclusions
homomorphisms
$$
\tilde H^\cdot(M^k)\to \tilde H^\cdot(\Sigma_{k,k+1})\to\tilde H
^\cdot(\Sigma_{1,k,k+1})
$$
is zero. We obtain that the homomorphism $\tilde H^\cdot(M^k)\to
\tilde H^\cdot(\Sigma_{k,k+1})$ is zero. \qquad $\square$

Let $X\subset\mathbb S$ be an open subset of $\mathbb S$ whose
compliment is a neighborhood deformation retract, we denote by
$\hat\delta:\tilde H^i(X)\to H^{i+1}(\mathbb B,X)$ natural
isomorphism induced by the exact sequence of the pair $\mathbb
B,X$.

Now consider the Mayer--Vietoris exact sequence of the pair
$\Sigma_{k,k+1},\,M^{k-1}$:
\begin{multline*}
\ldots H^{i-1}(M^k)\to \\ H^{i-1}(\Sigma_{k,k+1})\oplus
H^{i-1}(M^{k-1})\stackrel{\theta}{\to}H^{i-1}(\Sigma_{1,k,k+1})\stackrel{d}{\to}
H^i(M^k)\ldots
\end{multline*}
and its relative version:
\begin{multline*}
\ldots H^{i}(\mathbb B,M^k)\to \\ H^{i}(\mathbb
B,\Sigma_{k,k+1})\oplus H^{i}(\mathbb
B,M^{k-1})\stackrel{\theta}{\to}H^{i}(\mathbb B,
\Sigma_{1,k,k+1})\stackrel{d}{\to}H^{i+1}(\mathbb B,M^k)\ldots.
\end{multline*}
Then $\hat\delta$ establishes the isomorphism of these two exact
sequences. Moreover, Lemma~4 implies that long exact sequences
split in the short ones:
$$
0\to H^{i-1}(\Sigma_{k,k+1})\oplus
H^{i-1}(M^{k-1})\stackrel{\theta}{\to}H^{i-1}(\Sigma_{1,k,k+1})\stackrel{d}{\to}
H^i(M^k)\to 0 \eqno (3)
$$
and similarly for the relative version.

\begin{lemma} Let $\xi\in H^\cdot(\Sigma_{k,k+1}),\ \eta\in
H^\cdot(M^{k-1})$. Then $\hat\delta\xi\smile\hat\delta\eta=0$ if
and only if $\left(\xi\bigr|_{\Sigma_{1,k,k+1}}\smile
\eta\bigr|_{\Sigma_{1,k,k+1}}\right)\in\mathrm{im}\,\theta$.
\end{lemma}
{\bf Proof.} The Proposition from the Appendix~A implies:
$$
\hat\delta\xi\smile\hat\delta\gamma=\delta\circ
d\left(\xi\bigr|_{\Sigma_{1,k,k+1}}\smile\eta\bigr|_{\Sigma_{1,k,k+1}}\right).
$$
Now the statement of the Lemma follows from the fact that
$\hat\delta$ is an isomorphism and the sequence (3) is exact.
\qquad $\square$

The next step is to find $\mathrm{im}\,\theta$. Given $\xi\in
H^\cdot(\Sigma_{k,k+1}),\ \eta\in H^\cdot(M^{k-1})$, we have:
$$
\theta(\xi\oplus\eta)=\xi\bigr|_{\Sigma_{1,k,k+1}}-\eta\bigr|_{\Sigma_{1,k,k+1}}.
$$
According to Corollary~2, the restriction
$H^\cdot(\Sigma_{k,k+1})\to
H^\cdot(\Sigma_{k,k+1})\bigr|_{\Sigma_{1,k,k+1}}$ is injective and
$$
H^\cdot(\Sigma_{1,k,k+1})=H^\cdot(\Sigma_{k,k+1})\bigr|_{\Sigma_{1,k,k+1}}\oplus
\left(v_{k-1}\smile
H^\cdot(\Sigma_{k,k+1})\bigr|_{\Sigma_{1,k,k+1}}\right).
$$

Recall that $M^{k-1}$ has the homotopy type of the Thom space of a
vector bundle over $G_{k-1}(e)\subset M^{k-1}$ with the Thom class
$u_{k-1}\in H^{\nu_{k-1}}(M^{k-1})$. We consider the map
$\varrho_k:H^\cdot(G_{k-1}(e))\to H^\cdot(\Sigma_{k,k+1})$, where
$$
v_{k-1}\smile\varrho_k(\zeta)\bigr|_{\Sigma_{1,k,k+1}}=
\pi_v(u_{k-1}\smile\zeta)\bigr|_{\Sigma_{1,k,k+1}},\quad \forall
\zeta\in H^\cdot(G_{k-1}(e)). \eqno (4)
$$
The identity (4) uniquely defines $\varrho_k$. Moreover, the map
$\varrho_k$ is injective and
$$
\mathrm{im}\,\theta=H^\cdot(\Sigma_{k,k+1})\bigr|_{\Sigma_{1,k,k+1}}\oplus
\left(v_{k-1}\smile\mathrm{im}\,\varrho_k\bigr|_{\Sigma_{1,k,k+1}}\right).
\eqno (5)
$$
The space $\Sigma_{k,k+1}$ has the homotopy type of the Grassmannian
$Gr_k(n)$ while $G_{k-1}(e)$ is identified with the Grassmannian
$\{F\in Gr_{k-1}(n): F\subset e^\perp\}=Gr_{k-1}(n-1)$. We are going
to explicitly compute the map $\varrho_k$ in the bases provided by
the Schubert cells in the Grassmannians.

In what follows, we identify the manifold $$\Sigma_{k,k+1}\cap
M_{k-1}\cap W_k=
\left(\Sigma_{k,k+1}\setminus\Sigma_{1,k,k+1}\right)\cap W_k$$
with the Grassmannian $Gr_k(n)$, where $A\in\Sigma_{k,k+1}\cap
M_{k-1}\cap W_k$ is identified with the subspace $\{x\in\mathbb
R^n: Ax=\lambda_1(A)x\}$. Obviously, $\Sigma_{k,k+1}\cap
M_{k-1}\cap W_k$ is a homotopy retract of $\Sigma_{k,k+1}$. In
particular, $H^\cdot(\Sigma_{k,k+1})=H^\cdot(Gr_k(n))$.

Let $e_1=e,e_2,\ldots,e_n$ be an orthogonal basis of $\mathbb
R^n$. The closed Schubert cells in $Gr_k(n)$ associated to this
basis are cycles which give an additive basis of
$H_\cdot(Gr_k(n))$ (see Appendix B). We also consider the dual
{\it Schubert basis} of $H^\cdot(Gr_k(n))$. The Schubert cells of
dimension $r\ge 0$ are in the one-to-one correspondence with the
partitions of $r$ in no more than $k$ positive integral summands
in such a way that each summand does not exceed $n-k$.

Similarly, Schubert cells associated to the basis $e_2,\ldots,e_n$
of $e^\perp=\mathbb R^{n-1}$ give the Schubert basis of
$H^\cdot(Gr_{k-1}(n-1))=H^\cdot(G_{k-1}(e))$. The elements of
dimension $r$ of this basis are in the one-to-one correspondence
with the partitions of $r$ in less than $k$ positive integral
summands in such a way that each summand does not exceed $n-k$.

\begin{lemma} The map $\varrho_k:H^\cdot(G_{k-1}(e))\to
H^\cdot(Gr_k(n))$ sends the element of the Schubert basis of
$H^\cdot(G_{k-1}(e))$ associated to a partition in less than $k$
summands to the element of the Schubert basis of
$H^\cdot(Gr_k(n))$ associated to the same partition!
\end{lemma}
{\bf Proof.} We'll study the adjoint map
$\varrho^*_k:H_\cdot(Gr_k(n))\to H_\cdot(G_{k-1}(e))$. We have to
prove that $\varrho^*_k$ sends to zero the Schubert classes for
$H_\cdot(Gr_k(n))$ associated to the partitions in exactly $k$
summands, while the classes associated to the partitions in less
than $k$ summands are sent to the Schubert classes for
$H_\cdot(G_{k-1}(e))$ associated to the same partitions.

Let $C\subset Gr_k(n)$ be a Schubert cycle and $[C]$ its homology
class. We set $S^{\nu_{k-1}}_C=\bigcup\limits_{E\in
C}S^{\nu_{k-1}}_E$, c.f. (2). Then $\varrho^*_k[C]$ is the
homology class of the intersection of
$\phi_{k-1}\left(S^{\nu_{k-1}}_C\right)$ with
$G_{k-1}(e)=Gr_{k-1}(n-1)$. In other words, the map $\varrho^*_k$
is essentially determined by the set-valued map $\mathfrak
r_k:Gr_k(n)\dashrightarrow G_{k-1}(e)$, where $\mathfrak
r_k(E)=\phi_{k-1}\left(S^{\nu_{k-1}}_E\right)\cap G_{k-1}(e),\
E\in Gr_k(n)$.

It is easy to see that $\mathfrak r_k(E)=\{F\in Gr_{k-1}(n-1):
F\subset E\cap e^\perp\}.$ In particular, $\mathfrak r_k$ is
one-valued on $\{E\in Gr_k(n): E\not\perp e\}$; if $E\not\perp e$,
then $\mathfrak r_k(E)=E\cap e^\perp$. The (one-valued) map
$F\mapsto (F+\mathbb Re),\ F\in Gr_{k-1}(n-1),$ is a right inverse
of $\mathfrak r_k$.

Let $\mathfrak d$ be a starting from the unit Schubert symbol for
$Gr_k(n)$. Then $\mathfrak r_k\left(Sc^{\mathfrak
d}_k(n)\right)=Sc^{\mathfrak d'}_{k-1}(n-1)$, where $\mathfrak d'$
is obtained from $\mathfrak d$ by removing the first unit. Indeed,
$e\in E,\ \forall\,E\in Sc^{\mathfrak d}_k(n)$, and the desired
equality easily follows from the definitions.

Now assume that $\mathfrak d$ is a starting from 0 Schubert symbol
for $Gr_k(n)$. We'll show that $\mathfrak r_k\left(Sc^{\mathfrak
d}_k(n)\right)$ is contained in the union of Schubert cells whose
dimension is smaller than the dimension of $Sc^{\mathfrak
d}_k(n)$. This fact completes the proof of Lemma~6.

Let $F\in Sc^{\mathfrak d}_k(n)$ and $\hat F\in\mathfrak
r_k\left(Sc^{\hat{\mathfrak d}}_{k-1}(n-1)\right)$. Recall that
$$
d_i^{{\mathfrak d}}=\min\{j:\dim(E_j\cap F)=i\},\quad
i=1,\ldots,k,
$$
$$
d_i^{\hat{\mathfrak d}}=\min\{j:\dim(E_{j+1}\cap \hat F)=i\},\quad
i=1,\ldots,k-1,
$$
where $E_j=span\{e_1,\ldots,e_j\},\
j=1,\ldots,n.$ On the other hand,
$$
\dim(E_j\cap F_j)-1\le\dim(E_j\cap \hat F_j)\le\dim(E_j\cap F_j);
$$
hence $d_i^{{\mathfrak d}}\le d_i^{\hat{\mathfrak d}}+1\le
d_{i+1}^{{\mathfrak d}}$. Moreover, $d_1^{{\mathfrak d}}>1$, since
the symbol $\mathfrak d$ starts from 0. We obtain:
\begin{multline*}
\dim Sc^{\hat{\mathfrak d}}_{k-1}(n-1)= \\
\sum_{i=1}^{k-1}( d_i^{\hat{\mathfrak d}}-i)=
\sum_{i=1}^{k-1}\left(( d_i^{\hat{\mathfrak
d}}+1)-(i+1)\right)\le\sum_{i=2}^k(d_i^{{\mathfrak d}}-i)<\dim
Sc_k^{\mathfrak d}(n).  \end{multline*} $\square$

\begin{lemma} Let $w\in H^1(Gr_k(n))$ be the first
Stiefel--Whitney class of the tautological bundle and
$\xi\in\mathrm{im}\varrho_k$. Then
$(w\smile\xi)\in\mathrm{im}\,\varrho_k$ if and only if $\xi$ is a
sum of Schubert classes associated to partitions in less than
$k-1$ summands.
\end{lemma}
{\bf Proof.} Let $\Pi_j$ be the linear hull of the Schubert
classes associated to the partitions in exactly $j$ summands,
$j=0,1,\ldots,k,$ and $\xi$ the Schubert class associated to the
partition $a_1+\cdots+a_j$. The Pieri formula (see Appendix~B)
implies that the difference of $w\smile\xi$ and the Schubert class
associated to the partition $1+a_1+\cdots+a_j$ belongs to $\Pi_j$.
On the other hand, according to Lemma~6,
$\mathrm{im}\varrho_k=\bigoplus\limits_{j=0}^{k-1}\Pi_j.\qquad
\square$

\medskip Now we are ready to compute $\ker\mathbf d_k$ and thus
complete the proof of Theorem~2. Let $\xi\in H^\cdot(M^{k-1})$; then
$\xi=\hat\delta(u_{k-1}\smile\zeta)$ for a unique $\zeta\in
H^\cdot(G_{k-1}(e))$. We have:
$$
\mathbf
d_k(\xi)=\hat\delta(\gamma_k)\smile\hat\delta(u_{k-1}\smile\zeta),
$$
where $\gamma_k\in H^1(\Sigma_{k,k+1})$ was defined in Section~2.
According to Lemma~5, $\mathbf d_k(\xi)=0$ if and only if
$$
(u_{k-1}\smile\zeta)\bigr|_{\Sigma_{1,k,k+1}}\smile\gamma_k\bigr|_{\Sigma_{1,k,k+1}}=
\left(v_{k-1}\smile(\varrho_k(\zeta)\smile\gamma_k)\bigr|_{\Sigma_{1,k,k+1}}\right)\in
\mathrm{im}\,\theta.
$$
Further, $Gr_k(n)=\Sigma_{k,k+1}\cap M_{k-1}\cap W_k$ is a
homotopy retract of $\Sigma_{k,k+1}$, and
$\gamma_k\bigr|_{Gr_k(n)}$ is the first Stiefel--Whitney class of
the tautological bundle of the Grassmannian $Gr_k(n)$, i.\,e.
$\gamma_k\bigr|_{Gr_k(n)}=w$ (see Proposition~2). Now the equality
(5) implies that $\mathbf d_k(\xi)=0$ if and only if
$$
(\varrho_k(\zeta)\smile w)\in\mathrm{im}\varrho_k.
$$

It follows from Lemma~7 and the injectivity of $\varrho_k$ that
$\dim\ker\mathbf d_k$ is equal to the number of partitions in no
more than $k-2$ natural summands in such a way that each summand
does not exceed $n-k$. In other words, $\dim\ker\mathbf
d_k=\dbinom{n-2}{k-2}$. At the same time, the isomorphisms $
H^\cdot(\mathbb B,M^{k-1})\cong\tilde H^\cdot(M^{k-1})\cong
H^\cdot(Gr_{k-1}(n-1))$ imply that $\dim H^\cdot(\mathbb
B,M^{k-1})=\dbinom{n-1}{k-1}.$ The Pascal triangle identity
$\dbinom{n-1}{k-1}=\dbinom{n-2}{k-2}+\dbinom{n-2}{k-1}$ gives:
$$
\dim H^\cdot(\mathcal B,M^{k-1})=\dim\ker\mathbf
d_k+\dim\ker\mathbf d_{k+1}. \eqno \square
$$

\appendix

\section{A property of the cohomological product}

Let $M$ be a simplicial  complex and $X\subset M$ its subcomplex;
we denote by \linebreak $\delta_X: H^*(X)\to H^{*+1}(M,X)$ the
connecting homomorphism in the cohomological exact sequence of the
pair $M,X$. Let $Y\subset M$ be one more subcomplex and
$d:H^*(X\cap Y)\to H^{*+1}(X\cup Y)$ the connecting homomorphism
in cohomological Mayer--Vietoris exact sequence of the pair $X,Y$.

\begin{proposition} Let $\xi\in H^\cdot(X),\ \eta\in H^\cdot(Y)$; then
$$\delta_X\xi\smile\delta_Y\eta=\delta_{X\cup Y}\circ d\left(\xi|_{X\cap
Y}\smile\eta|_{X\cap Y}\right).$$
\end{proposition}
{\bf Proof.} We set $\zeta=d(\xi|_{X\cap Y}\smile\eta|_{X\cap
Y})$. Let $x$ and $y$ be cocycles representing cohomology classes
$\xi$ and $\eta$. Any cocycle $z$ such that
$$
z|_X=\delta u,\ z|_Y=\delta v,\quad u|_{X\cap Y}-v|_{X\cap
Y}=x|_{X\cap Y}\smile y|_{X\cap Y} \eqno (A)
$$
for some cochains $u,v$ is a representive of $\zeta$. We do as
follows: extend $x$ and $y$ to cochain $\hat x$ and $\hat y$
defined on $X\cup Y$ and set $z=\hat x\smile\delta\hat y$. Then
conditions (A) are satisfied for $u=(-1)^{\dim x}x\smile\hat
y|_X,\ v=0$ and we have: $\delta z=\delta\hat x\smile\delta\hat
y$. \qquad $\square$

\section{Schubert cells}

Schubert cells give cell complex structures of Grassmannians. They
are indexed by {\it Schubert symbols}. A Schubert symbol
$\mathfrak d$ for $Gr_k(n)$ is a sequence of zeros and units that
contains exactly $k$ units and $n-k$ zeros. The total number of
symbols for $Gr_k(n)$ (i.\,e. the number of cells in the cell
complex) is $\binom nk$. We denote by $d_i^{\mathfrak d}$ the
number of $i$th unit in the sequence; then $1\le d^{\mathfrak
d}_1<\cdots<d^{\mathfrak d}_k\le n$.

We treat simultaneously the real and complex cases. Let
$e_1,\ldots,e_n$ be a basis of $\mathbb R^n$ in the real case and
a basis of $\mathbb C^n$ in the complex case. We set
$$
E_i=span\{e_1,\ldots,e_i\},\quad i=1,\ldots,n.
$$
The Schubert cell $Sc_k^{\mathfrak d}(n)$ is defined as follows:
$$
Sc_k^{\mathfrak d}(n)=\left\{F\in Gr_k(n): \dim(F\cap
E_{d_i^{\mathfrak d}})=i,\ \dim(F\cap E_{d_i^{\mathfrak
d}-1})=i-1\right\}.
$$

There is a one-to-one correspondence between Schubert symbols for
$Gr_k(n)$ and partitions of nonnegative integers in no more than
$k$ positive integral summands in such a way that each summand
does not exceed $(n-k)$. The summands associated to the symbol
$\mathfrak d$ are numbers of zeros to the left of each unit
presented in the symbol. In other words, the summands are nonzero
terms of the sequence $(d^{\mathfrak d}_i-i),\ i=1,\ldots,k.$

The dimension of the Schubert cell associated to a partition of
the number $r$ is equal to $r$ in the real case and to $2r$ in the
complex case. We thus have:
$$
\dim Sc_k^{\mathfrak d}(n)=\epsilon\sum_{i=1}^k(d^{\frak d}_i-i),
$$
where $\epsilon=1$ in the real case and $\epsilon=2$ in the
complex case. The closure $\overline{Sc_k^{\mathfrak d}(n)}$ is a
cycle over $\mathbb Z_2$ in the real case and a cycle over
$\mathbb Z$ in the complex case (a {\it Schubert cycle}). In both
cases, the homology classes of the Schubert cycles form an
additive basis of the total homology group of $Gr_k(n)$ (over
$\mathbb Z_2$ in the real case and over $\mathbb Z$ in the complex
one). Moreover, in the complex case the homology groups are free.

The dual basis of the total cohomology group of $Gr_k(n)$ is
called the {\it Schubert basis}. Repeat that the elements of
dimension $r$ of this basis, the $r$-dimensional {\it Schubert
classes}, are in the one-to-one correspondence with partitions of
$r$ in no more than $k$ natural summands, where each summand does
not exceed $n-k$.

The Stiefel--Whitney (in the real case) and Chern (in the complex
case) characteristic classes of the tautological bundle are
Schubert classes associated to partitions in units. In particular,
the Stiefel--Whitney class $w_1$ in the real case and the Chern
class $c_1$ in the complex case are associated to the unique
``partition" of 1.

There is a useful {\it Pieri formula} which computes the
cohomological product of the Schubert class associated to a
partition with one summand $a$ and the Schubert class associated
to any partition $b_1+\cdots+b_j$, where $b_1\le\cdots\le b_j$.
The product equals the sum of all Schubert classes of dimension
$a+\sum\limits_{i=1}^jb_i$ associated to the partitions
$c_0+c_1+\cdots+c_j$ such that $b_{i-1}\le c_i\le b_i,\
i=1,\ldots,j,\ b_0=0$.

See details in \cite[Ch.5]{MiSt} and \cite[Ch.1.5]{GrHa}.

\end{document}